\documentclass[11pt, twoside]{article}
\usepackage{latexsym}
\usepackage{amsmath}
\usepackage{amssymb}
\usepackage[all]{xy}
\usepackage{amsfonts}
\usepackage{verbatim}
\usepackage{amsthm}
\usepackage{mathrsfs}
\usepackage{epsfig}
\usepackage{xy}
\usepackage{array}
\usepackage{stmaryrd}
\usepackage{graphicx,color}
\usepackage{xcolor}
\usepackage{tikz}
\usetikzlibrary{arrows,calc}
\usepackage{etex}
\usepackage{mathdots}
\usepackage{float}
\usepackage{graphics}
\usepackage{pdflscape}
\usepackage{extarrows}
\usepackage{anysize,hyperref}
\input xypic
\xyoption{all}
\usepackage{bm}
\usepackage[perpage,symbol]{footmisc}
\usepackage{setspace}
\def\A{\mathcal{A}}

\def\P{\mathcal{P}}
\def\I{\mathcal{I}}
\def\B{\mathcal{B}}
\def\C{\mathscr{C}}
\def\E{\mathbb{E}}

\def\s{\mathfrak{s}}
\def\id{\mathrm{id}}
\def\op{^\mathrm{op}}
\def\Ab{\mathsf{Ab}}
\def\del{\delta}
\def\dr{\ar@{->}[r]}

\def\X{\mathscr{X}}

\def\Hom{\mbox{Hom}}

\newcommand{\CC}{{\bf{C}}^{n+2}_{\C}}
\newcommand{\mr}{\hbox{\boldmath$\cdot$}}
\newcommand{\ov}{\overset}
\newcommand{\lra}{\longrightarrow}
\newcommand{\co}{\colon}
\newcommand{\uas}{^{\ast}}            
\newcommand{\sas}{_{\ast}}
\newcommand{\Xd}{\langle X^{\mr},\del\rangle}  
\newcommand{\Yr}{\langle Y^{\mr},\rho\rangle}  

\newcommand{\ush}{^\sharp}           
\newcommand{\ssh}{_\sharp}
\usepackage{enumitem}
\SelectTips{cm}{10}

\begin{document}
\baselineskip=15pt
\title{\Large{\bf From $\bm{n}$-exangulated categories to $\bm{n}$-abelian categories\footnotetext{Yu Liu was supported by the Fundamental Research Funds for the Central Universities (Grant No. 2682019CX51) and the National Natural Science Foundation of China (Grant No. 11901479). Panyue Zhou was supported by the National Natural Science Foundation of China (Grant Nos. 11901190 and 11671221), and the Hunan Provincial Natural Science Foundation of China (Grant No. 2018JJ3205),  and by the Scientific Research Fund of Hunan Provincial Education Department (Grant No. 19B239).}}}
\bigskip
\author{Yu Liu and Panyue Zhou}

\date{}

\maketitle
\def\blue{\color{blue}}
\def\red{\color{red}}

\newtheorem{theorem}{Theorem}[section]
\newtheorem{lemma}[theorem]{Lemma}
\newtheorem{corollary}[theorem]{Corollary}
\newtheorem{proposition}[theorem]{Proposition}
\newtheorem{conjecture}{Conjecture}
\theoremstyle{definition}
\newtheorem{definition}[theorem]{Definition}
\newtheorem{question}[theorem]{Question}
\newtheorem{remark}[theorem]{Remark}
\newtheorem{remark*}[]{Remark}
\newtheorem{example}[theorem]{Example}
\newtheorem{example*}[]{Example}
\newtheorem{condition}[theorem]{Condition}
\newtheorem{condition*}[]{Condition}
\newtheorem{construction}[theorem]{Construction}
\newtheorem{construction*}[]{Construction}

\newtheorem{assumption}[theorem]{Assumption}
\newtheorem{assumption*}[]{Assumption}

\baselineskip=17pt
\parindent=0.5cm

\begin{abstract}
\baselineskip=16pt
Herschend-Liu-Nakaoka introduced the notion of $n$-exangulated categories. It is not only a higher
dimensional analogue of extriangulated categories defined by Nakaoka-Palu,
but also gives a  simultaneous generalization of $n$-exact categories in the sense of
Jasso and $(n+2)$-angulated in the sense of Geiss-Keller-Oppermann.
Let $\C$ be an $n$-exangulated category with enough projectives and enough injectives, and $\X$ a cluster tilting subcategory of $\C$. In this article,  we show that the quotient category $\C/\X$ is an $n$-abelian category.
This extends a result of Zhou-Zhu for $(n+2)$-angulated
categories. Moreover, it highlights new phenomena when it is applied to $n$-exact categories.\\[0.5cm]
\medskip
\end{abstract}

\pagestyle{myheadings}
\markboth{\rightline {\scriptsize   Y. Liu and P. Zhou}}
         {\leftline{\scriptsize  From $n$-exangulated categories to $n$-abelian categories}}

\section{Introduction}
Cluster tilting theory permitted to construct abelian categories from some
triangulated categories. By Buan-Marsh-Reiten \cite[Theorem 2.2]{BMR} in cluster categories, by Keller-Reiten
\cite[Proposition 2.1]{KR} in the $2$-Calabi-Yau case, then by Koenig-Zhu \cite[Theorem 3.3]{KZ} and Iyama-Yoshino \cite[Corollary 6.5]{IY} in the general case, one can pass from triangulated categories to abelian categories by factoring out cluster tilting subcategories. Cluster tilting theory is also permitted to construct abelian categories from some
exact categories.  Demonet-Liu \cite{DL} provided a general framework for passing from exact
categories to abelian categories by factoring out cluster tilting subcategories. More precisely,
Demonet-Liu \cite{DL} showed that if  an exact category $\B$ with enough projectives and injectives has a cluster tilting
subcategory $\X$, then $\B/\X$ is an abelian category.

The notion of an extriangulated category was introduced in \cite{NP}, which is a simultaneous generalization of
exact category and triangulated category. In particular, exact
categories and triangulated categories are extriangulated categories.
Hence, many results hold on exact categories and triangulated categories can be unified in the same framework.
Of course, there are a lot of examples of extriangulated categories which are neither exact categories nor triangulated categories, see \cite{NP,ZZ1}.
Zhou and Zhu \cite{ZZ2} obtained abelian categories as quotients from triangulated categories modulo cluster tilting subcategories.
This unifies a result by Koenig-Zhu for triangulated categories and a result by Demonet-Liu for exact categories.

In \cite{GKO}, Geiss, Keller and Oppermann introduced $(n+2)$-angulated categories. There
are a higher dimensional analogue of triangulated categories,
 in the sense that triangles are replaced by $(n+2)$-angles, that is, morphism sequences of length $(n+2)$. 
 Thus a $1$-angulated category is precisely a triangulated category. 
 A typical example is a certain $n$-cluster tilting subcategory of a triangulated category.
 The notions of $n$-abelian and $n$-exact categories were introduced by Jasso in \cite{Ja}.
 An important source of examples of $n$-abelian and $n$-exact categories are $n$-cluster tilting subcategories.
 For $n=1$ they specialise to abelian and exact categories.
Zhou and Zhu \cite{ZZ3} proved that any quotient of an $(n+2)$-angulated
category modulo a cluster tilting subcategory is an $n$-abelian category.
When $n=1$, it is just the Koenig-Zhu's result.

Recently, Herschend-Liu-Nakaoka \cite{HLN} defined $n$-exangulated categories as a higher dimensional analogue
of extriangulated categories.  In particular, $n$-exangulated categories
simultaneously generalize $(n+2)$-angulated and $n$-exact categories.
They also gave some examples of $n$-exangulated categories, which are neither $n$-exact nor $(n+2)$-angulated.
For more examples of $n$-exangulated categories, see \cite{LZ}.

It is very natural to ask if any quotient of an $n$-exangulated
category modulo a cluster tilting subcategory is an $n$-abelian category.
In this article, we give  an affirmative answer.
We define cluster tilting subcategories of $n$-exangulated categories and prove the following main result.

\begin{theorem}{\rm (see Theorem \ref{main} for details)}
Let $\C$ be an $n$-exangulated category with enough projectives and enough injectives, and $\X$ a cluster tilting subcategory of $\C$.
Then $\C/\X$ is an $n$-abelian category.
\end{theorem}

Since an $(n+2)$-angulated category can be viewed
as an $n$-exangulated category with enough projectives and enough injectives,
when $\C$ is an $(n+2)$-angulated category, our main result is just the Theorem 3.4 in \cite{ZZ3}.
Moreover, this result is completely new when it is applied to $n$-exact categories.
At the same time, it also provides a new source to construct $n$-abelian categories.

This article is organized as follows. In Section 2, we review some elementary definitions
and facts about $n$-abelian, right $(n+2)$-angulated  and $n$-exangulated categories. In Section 3, we prove our main
result and apply to $(n+2)$-angulated and $n$-exact categories.

\section{Preliminaries}

In this section, we review basic concepts and results concerning $n$-abelian categoires, right $(n+2)$-angulated categories
and $n$-exangulated categories.

\subsection{$n$-abelian categories}
Let $\A$ be an additive category and $f\colon A\rightarrow B$ a morphism in $\A$. A \emph{weak cokernel} of $f$ is a morphism
$g\colon B\rightarrow C$ such that for any $X\in\A$ the sequence of abelian groups
$$\A(C,X)\xrightarrow{~\A(g,X)~}\A(B,X)\xrightarrow{~\A(X,f)~}\A(A,X)$$
is exact. Equivalently, $g$ is a weak cokernel of $f$ if $gf=0$ and for each morphism
$h\colon B\rightarrow X$ such that $hf=0$ there exists a (not necessarily unique) morphism
$p\colon C\rightarrow X$ such that $h=pg$. These properties are subsumed in the following
commutative diagram
$$\xymatrix{A\ar[r]^{f}\ar[dr]_{0} & B\ar[r]^{g}\ar[d]^{h}&C\ar@{-->}[dl]^{p}  \\
& X &}$$
Clearly, a weak cokernel $g$ of $f$ is a cokernel of $f$ if and only if $g$ is an epimorphism.
The concept of a \emph{weak kernel} is defined dually.

\begin{definition}\cite[Definiton 2.2]{Ja} and \cite[Definition 2.4]{L}
Let $\A$ be an additive category and $f_0\colon A_0\rightarrow A_1$  a morphism in
$\A$. An $n$-\emph{cokernel} of $f_0$ is a sequence
$$(f_1,f_2,\cdots,f_{n})\colon A_1\xrightarrow{~f_1~}A_2\xrightarrow{~f_2~}\cdots \xrightarrow
{~f_{n-1}~}A_n\xrightarrow
{~f_{n}~}A_{n+1}$$
such that the induced sequence of abelian groups
$$\xymatrix{0\xrightarrow{~~}\A(A_{n+1},B)\xrightarrow{~~} \A(A_{n},B)\xrightarrow{~~}\cdots \xrightarrow{~~}
\A(A_{1},B)\xrightarrow{~~}
\A(A_{0},B)}$$
is exact for each object $B\in\A$. That is, the morphism $f_i$ is a weak cokernel of $f_{i-1}$ for all $i=1,2,\cdots,n-1$ and $f_{n}$ is a cokernel of $f_{n-1}$. In this case, we say the sequence
$$A_0\xrightarrow{~f_0~}A_1\xrightarrow{~f_1~}A_2\xrightarrow{~f_2~}\cdots \xrightarrow
{~f_{n-1}~}A_n\xrightarrow
{~f_{n}~}A_{n+1} $$
is \emph{right $n$-exact}.

We can define \emph{$n$-kernel} and \emph{left $n$-exact} sequence dually. The sequence
$$A_0\xrightarrow{~f_0~}A_1\xrightarrow{~f_1~}A_2\xrightarrow{~f_2~}A_3\xrightarrow{~f_3~}\cdots\xrightarrow
{~f_{n-1}~}A_n\xrightarrow{~f_n~}A_{n+1}$$
is called \emph{$n$-exact} if it is
both right $n$-exact and left $n$-exact.
\end{definition}

\begin{definition}\cite[Definiton 3.1]{Ja}\label{def0}
Let $n$ be a positive integer. An \emph{$n$-abelian category} is an additive category $\A$
which satisfies the following axioms:
\begin{itemize}
\item[(A0)] The category $\A$ has split idempotents.

\item[(A1)] Every morphism in $\A$ has an $n$-kernel and an $n$-cokernel.

\item[(A2)] For every monomorphism $f_0\colon A_0\to A_1$ in $\A$ there exists an $n$-exact sequence:
$$A_0\xrightarrow{~f_0}A_1\xrightarrow{~f_1~}A_2
\xrightarrow{~f_2~}\cdots\xrightarrow{~f_3~}
A_{n-1}\xrightarrow{~f_{n-1}~}A_{n}
\xrightarrow{~f_n~}A_{n+1}.$$

\item[(A2)$^{\textrm{op}}$] For every epimorphism $g_n\colon B_n\to B_{n+1}$ in $\A$ there exists an $n$-exact sequence:
$$B_0\xrightarrow{~g_0}B_1\xrightarrow{~g_1~}B_2
\xrightarrow{~g_2~}\cdots\xrightarrow{~g_3~}
B_{n-1}\xrightarrow{~g_{n-1}~}B_{n}
\xrightarrow{~g_n~}B_{n+1}.$$
\end{itemize}
\end{definition}

\begin{remark} $1$-abelian categories are precisely abelian categories in the usual sense. It is easy
to see that abelian categories have split idempotents; thus, if $n=1$, then axiom (A0) in
Definition \ref{def0} is redundant.
\end{remark}

\subsection{Right $(n+2)$-angulated categories}

Let $\mathcal{A}$ be an additive category with an endofunctor $\Sigma^n:\mathcal{A}\rightarrow\mathcal{A}$. An $(n+2)$-$\Sigma^n$-$sequence$ in $\mathcal{A}$ is a sequence of morphisms
$$A_0\xrightarrow{f_0}A_1\xrightarrow{f_1}A_2\xrightarrow{f_2}\cdots\xrightarrow{f_{n-1}}A_n\xrightarrow{f_n}A_{n+1}\xrightarrow{f_{n+1}}\Sigma^n A_0.$$
Its {\em left rotation} is the $(n+2)$-$\Sigma^n$-sequence
$$A_1\xrightarrow{f_1}A_2\xrightarrow{f_2}A_3\xrightarrow{f_3}\cdots\xrightarrow{f_{n}}A_{n+1}\xrightarrow{f_{n+1}}\Sigma^n A_0\xrightarrow{(-1)^{n}\Sigma^n f_0}\Sigma^n A_1.$$
A \emph{morphism} of $(n+2)$-$\Sigma^n$-sequences is  a sequence of morphisms $\varphi=(\varphi_0,\varphi_1,\cdots,\varphi_{n+1})$ such that the following diagram commutes
$$\xymatrix{
A_0 \ar[r]^{f_0}\ar[d]^{\varphi_0} & A_1 \ar[r]^{f_1}\ar[d]^{\varphi_1} & A_2 \ar[r]^{f_2}\ar[d]^{\varphi_2} & \cdots \ar[r]^{f_{n}}& A_{n+1} \ar[r]^{f_{n+1}}\ar[d]^{\varphi_{n+1}} & \Sigma^n A_0 \ar[d]^{\Sigma^n \varphi_0}\\
B_0 \ar[r]^{g_0} & B_1 \ar[r]^{g_1} & B_2 \ar[r]^{g_2} & \cdots \ar[r]^{g_{n}}& B_{n+1} \ar[r]^{g_{n+1}}& \Sigma^n B_0
}$$
where each row is an $(n+2)$-$\Sigma^n$-sequence. It is an {\em isomorphism} if $\varphi_0, \varphi_1, \varphi_2, \cdots, \varphi_{n+1}$ are all isomorphisms in $\mathcal{A}$.
\medskip

We recall the notion of a right $(n+2)$-angulated category from \cite[Definition 2.1]{L}.
Compare with \cite[Definition 2.1]{L}, the condition (RN1)(a) is slightly different from that in \cite{L}, we don't assume that
the class $\Theta$ is closed under direct summands.

\begin{definition}\label{d1}
A {\em right} $(n+2)$-\emph{angulated category} is a triple $(\mathcal{A}, \Sigma^n, \Theta)$, where $\mathcal{A}$ is an additive category, $\Sigma^n$ is an endofunctor of $\mathcal{A}$ (which is called the $n$-suspension functor), and $\Theta$ is a class of $(n+2)$-$\Sigma^n$-sequences (whose elements are called right $(n+2)$-angles), which satisfies the following axioms:
\begin{itemize}
\item[(RN1)]
\begin{itemize}
\item[(a)] The class $\Theta$ is closed under isomorphisms and direct sums.

\item[(b)] For each object $A\in\mathcal{A}$ the trivial sequence
$$0\rightarrow A\xrightarrow{1_A}A\rightarrow 0\rightarrow\cdots\rightarrow 0\rightarrow 0$$
belongs to $\Theta$.

\item[(c)] Each morphism $f_0\colon A_0\rightarrow A_1$ in $\A$ can be extended to a right $(n+2)$-angle: $$A_0\xrightarrow{f_0}A_1\xrightarrow{f_1}A_2\xrightarrow{f_2}\cdots\xrightarrow{f_{n-1}}A_n\xrightarrow{f_n}A_{n+1}\xrightarrow{f_{n+1}}\Sigma^n A_0.$$
\end{itemize}
\item[(RN2)] If an $(n+2)$-$\Sigma^n$-sequence belongs to $\Theta$, then its left rotation belongs to $\Theta$.

\item[(RN3)] Each solid commutative diagram
$$\xymatrix{
A_0 \ar[r]^{f_0}\ar[d]^{\varphi_0} & A_1 \ar[r]^{f_1}\ar[d]^{\varphi_1} & A_2 \ar[r]^{f_2}\ar@{-->}[d]^{\varphi_2} & \cdots \ar[r]^{f_{n}}& A_{n+1} \ar[r]^{f_{n+1}}\ar@{-->}[d]^{\varphi_{n+1}} & \Sigma^n A_0 \ar[d]^{\Sigma^n \varphi_0}\\
B_0 \ar[r]^{g_0} & B_1 \ar[r]^{g_1} & B_2 \ar[r]^{g_2} & \cdots \ar[r]^{g_{n}}& B_{n+1} \ar[r]^{g_{n+1}}& \Sigma^n B_0
}$$ with rows in $\Theta$ can be completed to a morphism of  $(n+2)$-$\Sigma^n$-sequences.

\item[(RN4)] Given a commutative diagram
$$\xymatrix{
A_0\ar[r]^{f_0}\ar@{=}[d] & A_1 \ar[r]^{f_1}\ar[d]^{\varphi_1} & A_2 \ar[r]^{f_2} & \cdots\ar[r]^{f_{n-1}} & A_{n}\ar[r]^{f_{n}\quad} & A_{n+1} \ar[r]^{f_{n+1}} & \Sigma^n A_0\ar@{=}[d] \\
A_0\ar[r]^{g_0} & B_1 \ar[r]^{g_1}\ar[d]^{h_1} & B_2\ar[r]^{g_2} & \cdots\ar[r]^{g_{n-1}} & B_{n}\ar[r]^{g_{n}\quad} & B_{n+1} \ar[r]^{g_{n+1}} & \Sigma^n A_0\\
& C_2\ar[d]^{h_2} & & & & & \\
& \vdots\ar[d]^{h_{n-1}} & & & & & \\
& C_{n}\ar[d]^{h_{n}} & & & & & \\
& C_{n+1}\ar[d]^{h_{n+1}} & & & & & \\
& \Sigma^n A_1 & & & & & \\
}$$
whose top rows and second column belong to $\Theta$. Then there exist morphisms $\varphi_i\colon A_i\rightarrow B_i\ (i=2,3,\cdots,n+1)$, $\psi_j\colon B_j\rightarrow C_j\ (j=2,3,\cdots,n+1)$ and $\phi_k\colon A_k\rightarrow C_{k-1}\ (k=3,4,\cdots,n+1)$ with the following two properties:

(I) The sequence $(1_{A_1},\varphi_1, \varphi_2,\cdots,\varphi_{n+1})$ is a morphism of $(n+2)$-$\Sigma^n$-sequences.

(II) The $(n+2)$-$\Sigma^n$-sequence
$$A_2\xrightarrow{\left(
                    \begin{smallmatrix}
                      f_2 \\
                      \varphi_2 \\
                    \end{smallmatrix}
                  \right)} A_3\oplus B_2\xrightarrow{\left(
                             \begin{smallmatrix}
                               -f_3 & 0 \\
                               \varphi_3 & -g_2 \\
                               \phi_3 & \psi_2 \\
                             \end{smallmatrix}
                           \right)}
 A_4\oplus B_3\oplus C_2\xrightarrow{\left(
                                       \begin{smallmatrix}
                                         -f_4 & 0 & 0 \\
                                         -\varphi_4 & -g_3 & 0 \\
                                         \phi_4 & \psi_3 & h_2 \\
                                       \end{smallmatrix}
                                     \right)}A_5\oplus B_4\oplus C_3$$
$$\xrightarrow{\left(
                                       \begin{smallmatrix}
                                         -f_5 & 0 & 0 \\
                                         \varphi_5 & -g_4 & 0 \\
                                         \phi_5 & \psi_4 & h_3 \\
                                       \end{smallmatrix}
                                     \right)}\cdots\xrightarrow{\scriptsize\left(\begin{smallmatrix}
             -f_{n} & 0 & 0 \\
             (-1)^{n+1}\varphi_{n} & -g_{n-1} & 0 \\
             \phi_{n} & \psi_{n-1} & h_{n-2} \\
             \end{smallmatrix}
             \right)}A_{n+1}\oplus B_{n}\oplus C_{n-1}$$
$$\xrightarrow{\left(
                                                       \begin{smallmatrix}
                                                         (-1)^{n+1}\varphi_{n+1} &-g_{n} &0 \\
                                                          \phi_{n+1}& \psi_{n}& h_{n-1} \\
                                                       \end{smallmatrix}
                                                     \right)}B_{n+1}\oplus C_{n}\xrightarrow{(\psi_{n+1},\ h_{n})}C_{n+1}\xrightarrow{\Sigma^n f_1\circ h_{n+1}}\Sigma^n A_2 \hspace{10mm}$$
belongs to $\Theta$, and  $h_{n+1}\circ\psi_{n+1}=\Sigma^n f_0\circ g_{n+1}$.
   \end{itemize}
\end{definition}
The notion of a \emph{left $n$-angulated category} is defined dually.
\vspace{1mm}

If $\Sigma^n$ is an automorphism and the class $\Theta$ is closed under direct summands, it is easy to see that the converse of an axiom (RN2) also holds, thus the right
$n$-angulated category $(\A,\Sigma^n, \Theta)$ is an $n$-angulated category in the sense of Geiss-Keller-Oppermann \cite[Definition 1.1]{GKO} and in the sense of Bergh-Thaule \cite[Theorem 4.4]{BT}. If $(\A,\Sigma^n, \Theta)$ is a right $n$-angulated category, $(\A,\Omega^n, \Phi)$
is a left $n$-angulated category, $\Omega^n$ is a quasi-inverse of $\Sigma^n$ and $\Theta=\Phi$, then $(\A,\Sigma^n, \Theta)$ is an $n$-angulated category.

The following two lemmas can be found in \cite{ZZ3}, which will be used in the sequel.

\begin{lemma}\label{bu1}
Let $(\mathcal{A}, \Sigma^n, \Theta)$ be a right $(n+2)$-angulated category and
$$A_0\xrightarrow{f_0}A_1\xrightarrow{f_1}A_2\xrightarrow{f_2}\cdots\xrightarrow{f_{n-1}}A_n\xrightarrow{f_n}A_{n+1}\xrightarrow{f_{n+1}}\Sigma^n A_0.$$
be a right $(n+2)$-angle. Then
$f_{i}\circ f_{i-1}=0, ~i=1,2,\cdots,n+1$. That is to say, any composition of
two consecutive morphisms in a right $(n+2)$-angle vanishes.
\end{lemma}

\begin{lemma}\label{bu2}
Let $(\mathcal{A}, \Sigma^n, \Theta)$ be a right $(n+2)$-angulated category and
$$A_0\xrightarrow{f_0}A_1\xrightarrow{f_1}A_2\xrightarrow{f_2}\cdots\xrightarrow{f_{n-1}}A_n\xrightarrow{f_n}A_{n+1}\xrightarrow{f_{n+1}}\Sigma^n A_0.$$
be a right $(n+2)$-angle. Then $f_{i}$ is a weak cokernal of $f_{i-1}$, $i=1,2,\cdots,n+1$.
\end{lemma}

\subsection{$n$-exangulated categories}
Assume that $\C$ is an additive category.  Let us
recall definitions and properties related to $n$-exangulated categories from \cite{HLN}.

\begin{definition}\cite[Definition 2.1]{HLN}
Suppose that $\C$ is equipped with an additive bifunctor $\E\colon\C\op\times\C\to\Ab$, where $\Ab$ is the category of abelian groups. For any pair of objects $A,C\in\C$, an element $\del\in\E(C,A)$ is called an {\it $\E$-extension} or simply an {\it extension}. We also write such $\del$ as ${}_A\del_C$ when we indicate $A$ and $C$.

Let ${}_A\del_C$ be any extension. Since $\E$ is a bifunctor, for any $a\in\C(A,A')$ and $c\in\C(C',C)$, we have extensions
$$ \E(C,a)(\del)\in\E(C,A')\ \ \text{and}\ \ \E(c,A)(\del)\in\E(C',A). $$
We abbreviately denote them by $a_{\ast}\del$ and $c^{\ast}\del$.
In this terminology, we have
$$\E(c,a)(\del)=c^{\ast}a_{\ast}\del=a_{\ast}c^{\ast}\del\in\E(C',A').$$
For any $A,C\in\C$, the zero element ${}_A0_C=0\in\E(C,A)$ is called the {\it split $\E$-extension}.
\end{definition}

\begin{definition}\cite[Definition 2.3]{HLN}
Let ${}_A\del_C,{}_{A'}\del'_{C'}$ be any pair of $\E$-extensions. A {\it morphism} $(a,c)\colon\del\to\del'$ of extensions is a pair of morphisms $a\in\C(A,B)$ and $c\in\C(A',C')$ in $\C$, satisfying the equality
$$a_{\ast}\del=c^{\ast}\del'. $$
\end{definition}

\begin{definition}\cite[Definition 2.7]{HLN}
Let $\bf{C}_{\C}$ be the category of complexes in $\C$. As its full subcategory, define $\CC$ to be the category of complexes in $\C$ whose components are zero in the degrees outside of $\{0,1,\ldots,n+1\}$. Namely, an object in $\CC$ is a complex $X^{\mr}=\{X_i,d^X_i\}$ of the form
\[ X_0\xrightarrow{d^X_0}X_1\xrightarrow{d^X_1}\cdots\xrightarrow{d^X_{n-1}}X_n\xrightarrow{d^X_n}X_{n+1}. \]
We write a morphism $f^{\mr}\co X^{\mr}\to Y^{\mr}$ simply $f^{\mr}=(f_0,f_1,\ldots,f_{n+1})$, only indicating the terms of degrees $0,\ldots,n+1$.
\end{definition}

\begin{definition}\cite[Definition 2.11]{HLN}
By Yoneda lemma, any extension $\del\in\E(C,A)$ induces natural transformations
\[ \del\ssh\colon\C(-,C)\Rightarrow\E(-,A)\ \ \text{and}\ \ \del\ush\colon\C(A,-)\Rightarrow\E(C,-). \]
For any $X\in\C$, these $(\del\ssh)_X$ and $\del\ush_X$ are given as follows.
\begin{enumerate}
\item[\rm(1)] $(\del\ssh)_X\colon\C(X,C)\to\E(X,A)\ ;\ f\mapsto f\uas\del$.
\item[\rm (2)] $\del\ush_X\colon\C(A,X)\to\E(C,X)\ ;\ g\mapsto g\sas\delta$.
\end{enumerate}
We abbreviately denote $(\del\ssh)_X(f)$ and $\del\ush_X(g)$ by $\del\ssh(f)$ and $\del\ush(g)$, respectively.
\end{definition}

\begin{definition}\cite[Definition 2.9]{HLN}
 Let $\C,\E,n$ be as before. Define a category $\AE:=\AE^{n+2}_{(\C,\E)}$ as follows.
\begin{enumerate}
\item[\rm(1)]  A pair $\Xd$ is an object of the category $\AE$ with $X^{\mr}\in\CC$
and $\del\in\E(X_{n+1},X_0)$ is called an $\E$-attached
complex of length $n+2$, if it satisfies
$$(d_0^X)_{\ast}\del=0~~\textrm{and}~~(d^X_n)^{\ast}\del=0.$$
We also denote it by
$$X_0\xrightarrow{d_0^X}X_1\xrightarrow{d_1^X}\cdots\xrightarrow{d_{n-2}^X}X_{n-1}
\xrightarrow{d_{n-1}^X}X_n\xrightarrow{d_n^X}X_{n+1}\overset{\delta}{\dashrightarrow}$$
\item[\rm (2)]  For such pairs $\Xd$ and $\langle Y^{\mr},\rho\rangle$,  $f^{\mr}\colon\Xd\to\langle Y^{\mr},\rho\rangle$ is
defined to be a morphism in $\AE$ if it satisfies $(f_0)_{\ast}\del=(f_{n+1})^{\ast}\rho$.

\end{enumerate}
\end{definition}

\begin{definition}\cite[Definition 2.13]{HLN}
 An {\it $n$-exangle} is an object $\Xd$ in $\AE$ that satisfies the listed conditions.
\begin{enumerate}
\item[\rm (1)] The following sequence of functors $\C\op\to\Ab$ is exact.
$$
\C(-,X_0)\xLongrightarrow{\C(-,\ d^X_0)}\cdots\xLongrightarrow{\C(-,\ d^X_n)}\C(-,X_{n+1})\xLongrightarrow{~\del\ssh~}\E(-,X_0)
$$
\item[\rm (2)] The following sequence of functors $\C\to\Ab$ is exact.
$$
\C(X_{n+1},-)\xLongrightarrow{\C(d^X_n,\ -)}\cdots\xLongrightarrow{\C(d^X_0,\ -)}\C(X_0,-)\xLongrightarrow{~\del\ush~}\E(X_{n+1},-)
$$
\end{enumerate}
In particular any $n$-exangle is an object in $\AE$.
A {\it morphism of $n$-exangles} simply means a morphism in $\AE$. Thus $n$-exangles form a full subcategory of $\AE$.
\end{definition}

\begin{definition}\cite[Definition 2.22]{HLN}
Let $\s$ be a correspondence which associates a homotopic equivalence class $\s(\del)=[{}_AX^{\mr}_C]$ to each extension $\del={}_A\del_C$. Such $\s$ is called a {\it realization} of $\E$ if it satisfies the following condition for any $\s(\del)=[X^{\mr}]$ and any $\s(\rho)=[Y^{\mr}]$.
\begin{itemize}
\item[{\rm (R0)}] For any morphism of extensions $(a,c)\co\del\to\rho$, there exists a morphism $f^{\mr}\in\CC(X^{\mr},Y^{\mr})$ of the form $f^{\mr}=(a,f_1,\ldots,f_n,c)$. Such $f^{\mr}$ is called a {\it lift} of $(a,c)$.
\end{itemize}
In such a case, we simple say that \lq\lq$X^{\mr}$ realizes $\del$" whenever they satisfy $\s(\del)=[X^{\mr}]$.

Moreover, a realization $\s$ of $\E$ is said to be {\it exact} if it satisfies the following conditions.
\begin{itemize}
\item[{\rm (R1)}] For any $\s(\del)=[X^{\mr}]$, the pair $\Xd$ is an $n$-exangle.
\item[{\rm (R2)}] For any $A\in\C$, the zero element ${}_A0_0=0\in\E(0,A)$ satisfies
\[ \s({}_A0_0)=[A\ov{\id_A}{\lra}A\to0\to\cdots\to0\to0]. \]
Dually, $\s({}_00_A)=[0\to0\to\cdots\to0\to A\ov{\id_A}{\lra}A]$ holds for any $A\in\C$.
\end{itemize}
Note that the above condition {\rm (R1)} does not depend on representatives of the class $[X^{\mr}]$.
\end{definition}

\begin{definition}\cite[Definition 2.23]{HLN}
Let $\s$ be an exact realization of $\E$.
\begin{enumerate}
\item[\rm (1)] An $n$-exangle $\Xd$ is called a $\s$-{\it distinguished} $n$-exangle if it satisfies $\s(\del)=[X^{\mr}]$. We often simply say {\it distinguished $n$-exangle} when $\s$ is clear from the context.
\item[\rm (2)]  An object $X^{\mr}\in\CC$ is called an {\it $\s$-conflation} or simply a {\it conflation} if it realizes some extension $\del\in\E(X_{n+1},X_0)$.
\item[\rm (3)]  A morphism $f$ in $\C$ is called an {\it $\s$-inflation} or simply an {\it inflation} if it admits some conflation $X^{\mr}\in\CC$ satisfying $d_X^0=f$.
\item[\rm (4)]  A morphism $g$ in $\C$ is called an {\it $\s$-deflation} or simply a {\it deflation} if it admits some conflation $X^{\mr}\in\CC$ satisfying $d_X^n=g$.
\end{enumerate}
\end{definition}

\begin{definition}\cite[Definition 2.27]{HLN}
For a morphism $f^{\mr}\in\CC(X^{\mr},Y^{\mr})$ satisfying $f^0=\id_A$ for some $A=X_0=Y_0$, its {\it mapping cone} $M_f^{\mr}\in\CC$ is defined to be the complex
\[ X_1\xrightarrow{d^{M_f}_0}X_2\oplus Y_1\xrightarrow{d^{M_f}_1}X_3\oplus Y_2\xrightarrow{d^{M_f}_2}\cdots\xrightarrow{d^{M_f}_{n-1}}X_{n+1}\oplus Y_n\xrightarrow{d^{M_f}_n}Y_{n+1} \]
where $d^{M_f}_0=\begin{bmatrix}-d^X_1\\ f_1\end{bmatrix},$
$d^{M_f}_i=\begin{bmatrix}-d^X_{i+1}&0\\ f_{i+1}&d^Y_i\end{bmatrix}\ (1\le i\le n-1),$
$d^{M_f}_n=\begin{bmatrix}f_{n+1}&d^Y_n\end{bmatrix}$.

{\it The mapping cocone} is defined dually, for morphisms $h^{\mr}$ in $\CC$ satisfying $h_{n+1}=\id$.
\end{definition}

\begin{definition}\cite[Definition 2.32]{HLN}
An {\it $n$-exangulated category} is a triplet $(\C,\E,\s)$ of additive category $\C$, additive bifunctor $\E\co\C\op\times\C\to\Ab$, and its exact realization $\s$, satisfying the following conditions.
\begin{itemize}
\item[{\rm (EA1)}] Let $A\ov{f}{\lra}B\ov{g}{\lra}C$ be any sequence of morphisms in $\C$. If both $f$ and $g$ are inflations, then so is $g\circ f$. Dually, if $f$ and $g$ are deflations then so is $g\circ f$.

\item[{\rm (EA2)}] For $\rho\in\E(D,A)$ and $c\in\C(C,D)$, let ${}_A\langle X^{\mr},c\uas\rho\rangle_C$ and ${}_A\Yr_D$ be distinguished $n$-exangles. Then $(\id_A,c)$ has a {\it good lift} $f^{\mr}$, in the sense that its mapping cone gives a distinguished $n$-exangle $\langle M^{\mr}_f,(d^X_0)\sas\rho\rangle$.
\item[{\rm (EA2$\op$)}] Dual of {\rm (EA2)}.
\end{itemize}
Note that the case $n=1$, a triplet $(\C,\E,\s)$ is a  $1$-exangulated category if and only if it is an extriangulated category, see \cite[Proposition 4.3]{HLN}.
\end{definition}

\begin{example}
From \cite[Proposition 4.34]{HLN} and \cite[Proposition 4.5]{HLN},  we know that $n$-exact categories and $(n+2)$-angulated categories are $n$-exangulated categories.
There are some other examples of $n$-exangulated categories
 which are neither $n$-exact nor $(n+2)$-angulated, see \cite[Section 6]{HLN} for more details and \cite[Remark 4.5]{LZ}.
\end{example}

\begin{lemma}\label{a1}
Let $(\C,\E,\s)$ be an $n$-exangulated category, and
$$A_0\xrightarrow{\alpha_0}A_1\xrightarrow{\alpha_1}A_2\xrightarrow{\alpha_2}\cdots\xrightarrow{\alpha_{n-2}}A_{n-1}
\xrightarrow{\alpha_{n-1}}A_n\xrightarrow{\alpha_n}A_{n+1}\overset{\delta}{\dashrightarrow}$$
a distinguished $n$-exangle. Then we have the following long exact sequences:
$$\C(-, A_0)\xrightarrow{}\C(-, A_1)\xrightarrow{}\cdots\xrightarrow{}
\C(-, A_{n+1})\xrightarrow{}\E(-, A_{0})\xrightarrow{}\E(-, A_{1})\xrightarrow{}\E(-, A_{2});$$
$$\C(A_{n+1},-)\xrightarrow{}\C(A_{n},-)\xrightarrow{}\cdots\xrightarrow{}
\C(A_0,-)\xrightarrow{}\E(A_{n+1},-)\xrightarrow{}\E(A_{n},-)\xrightarrow{}\E(A_{n-1},-).$$
In particular, $\alpha_{i+1}$ is a weak cokernel of $\alpha_i$ and
$\alpha_{i}$ is a weak kernel of $\alpha_{i+1}$, for each $i\in\{1,2,\cdots,n\}$.
\end{lemma}

\proof This follows from Definition 2.13 and Corollary 3.11 in \cite{HLN}.  \qed

\begin{lemma}\label{factor}{\rm\cite[Lemma 3.3]{ZW}}
Let
$$\xymatrix{
A_0\ar[r]^{\alpha_0}\ar[d]^{f_0} &A_1 \ar[r]^{\alpha_2}\ar[d]^{f_1} &\cdot\cdot\cdot \ar[r]^{\alpha_{n-1}}&A_n\ar[d]^{f_n} \ar[r]^{\alpha_{n+1}} &A_{n+1} \ar[d]^{f_{n+1}} \ar@{-->}[r]^-{\delta} &\\
B_0\ar[r]^{\beta_0} &B_1\ar[r]^{\beta_2} &\cdot\cdot\cdot \ar[r]^{\beta_{n-1}} &B_n\ar[r]^{\beta_n}  &B_{n+1}\ar@{-->}[r]^-{\delta'} &
}
$$
be any morphism of distinguished $n$-exangles. Then
$f_0$ factors through $\alpha_0$
if and only if $(f_0)_{\ast}\delta=f_{n+1}^{\ast}\delta'=0$
if and only if $f_{n+1}$
factors through $\beta_n$.

In particular, in the case $\delta=\delta'$ and
$(f_0, f_1,\cdots,f_{n+1})=(1_{A_0},1_{A_1},\cdots,1_{A_{n+1}})$,
$\alpha_0$ is a section if and only if $\delta=0$ if and only if $\alpha_{n+1}$ is a retraction.
\end{lemma}

\begin{lemma}\label{inflation}
Let $(\C,\E,\s)$ be an $n$-exangulated category
and $\alpha\in\C(A,B), h\in\C(A,C)$ any pair of morphisms.
If $\alpha$ is an inflation, then so is $\left[
              \begin{smallmatrix}
                -\alpha\\ h
              \end{smallmatrix}
            \right]\in\C(A,B\oplus C)$.
\end{lemma}

\proof  Since $\alpha$ is an inflation, there exists a distinguished $n$-exangle
$$A\xrightarrow{\alpha}B\xrightarrow{g_1}A_2\xrightarrow{g_2}A_3
 \xrightarrow{g_3}\cdots\xrightarrow{g_{n-1}}A_n\xrightarrow{g_n}A_{n+1}\overset{\delta}{\dashrightarrow}$$
in $\C$. By (EA2$^{\rm op}$), we can observe that $(h,1_{A_{n+1}})$ has a good lift
$f^{\mr}=(h,h_1,h_2,\cdots,h_{n},1_{A_{n+1}})$, that is,  we have the commutative diagram
$$\xymatrix{
A\ar[r]^{\alpha}\ar[d]^h&B\ar[r]^{g_1}\ar[d]^{h_1} & A_2 \ar[r]^{g_2}\ar[d]^{h_2}  & \cdots \ar[r]^{g_{n-1}}& A_n \ar[r]^{g_n}\ar[d]^{h_n}&A_{n+1}\ar@{-->}[r]^{\;\;\delta}\ar@{=}[d] & \\
C\ar[r]^{f_0}&C_1 \ar[r]^{f_1} & C_2 \ar[r]^{f_2}  & \cdots \ar[r]^{f_{n-1}} & C_n \ar[r]^{f_n}&A_{n+1}\ar@{-->}[r]^{\;\;\eta}&}$$
of distinguished $n$-exangles, and the mapping cocone
$$A\xrightarrow{\left[
              \begin{smallmatrix}
                -\alpha\\ h
              \end{smallmatrix}
            \right]}B\oplus C\xrightarrow{\left[
              \begin{smallmatrix}
                -g_1&0\\ h_1&f_0
              \end{smallmatrix}
            \right]}A_2\oplus C_1\xrightarrow{\left[
              \begin{smallmatrix}
                -g_2&0\\ h_2&f_1
              \end{smallmatrix}
            \right]}\cdots\xrightarrow{\left[
              \begin{smallmatrix}
                -g_{n-1}&0\\ h_{n-1}&f_{n-2}
              \end{smallmatrix}
            \right]}A_{n}\oplus B_{n-1}\xrightarrow{\left[
              \begin{smallmatrix}
                h_n&f_{n-1}
              \end{smallmatrix}
            \right]}B_n\overset{f^{\ast}_n\delta}{\dashrightarrow}$$
is a distinguished $n$-exangle.  It follows that $\left[
              \begin{smallmatrix}
                -\alpha\\ h
              \end{smallmatrix}
            \right]\colon A\to B\oplus C$ is an inflation.  \qed

\begin{definition}\label{def2}\cite[Definition 3.14 and Definition 3.15]{ZW}
Let $(\C,\E,\s)$ be an $n$-exangulated category.
\begin{itemize}
\item[(1)] An object $P\in\C$ is called \emph{projective} if, for any distinguished $n$-exangle
$$A_0\xrightarrow{\alpha_0}A_1\xrightarrow{\alpha_1}A_2\xrightarrow{\alpha_2}\cdots\xrightarrow{\alpha_{n-2}}A_{n-1}
\xrightarrow{\alpha_{n-1}}A_n\xrightarrow{\alpha_n}A_{n+1}\overset{\delta}{\dashrightarrow}$$
and any morphism $c$ in $\C(P,A_{n+1})$, there exists a morphism $b\in\C(P,A_n)$ satisfying $\alpha_n\circ b=c$.
We denote the full subcategory of projective objects in $\C$ by $\P$.
Dually, the full subcategory of injective objects in $\C$ is denoted by $\I$.

\item[(2)] We say that $\C$ {\it has enough  projectives} if
for any object $C\in\C$, there exists a distinguished $n$-exangle
$$B\xrightarrow{\alpha_0}P_1\xrightarrow{\alpha_1}P_2\xrightarrow{\alpha_2}\cdots\xrightarrow{\alpha_{n-2}}P_{n-1}
\xrightarrow{\alpha_{n-1}}P_n\xrightarrow{\alpha_n}C\overset{\delta}{\dashrightarrow}$$
satisfying $P_1,P_2,\cdots,P_n\in\P$. We can define the notion of having \emph{enough injectives} dually.
\end{itemize}
\end{definition}

\begin{remark}
~\begin{itemize}
\item[\rm (1)]  In the case $n=1$, they agree with
the usual definitions \cite[Definition 3.23, Definition 3.25 and Definition 7.1]{NP}.

\item[\rm (2)] If $(\C,\E,\s)$ is an $n$-exact category, then they agree with
\cite[Definition 3.11, Definition 5.3 and Definition 5.5]{Ja}.

\item[\rm (3)] If $(\C,\E,\s)$ is an $(n+2)$-angulated category, then
$\P=\I$ consists of zero objects. Moreover it always has enough projectives and enough injectives.
\end{itemize}
\end{remark}

\begin{lemma}\label{projective}
Let $(\C,\E,\s)$ be an $n$-exangulated category. Then
an object $P\in\C$ is projective if and only if it satisfies
$\E(P,A)=0$ for any $A\in\C$.
\end{lemma}

\proof Assume that $P\in\C$ is a projective object.
For any $\delta\in\E(P,A)$, there exists a distinguished $n$-exangle
 $$A\xrightarrow{g_0}A_1\xrightarrow{g_1}A_2\xrightarrow{g_2}A_3
 \xrightarrow{g_3}\cdots\xrightarrow{g_{n-1}}A_n\xrightarrow{g_n}P\overset{\delta}{\dashrightarrow}$$
in $\C$. Since $P$ is a projective object, we have that $g_n$ is a retraction.
By Lemma \ref{factor}, we get that $\delta=0$.

Conversely, suppose $\E(P,A)=0$ for any $A\in\C$.
By Lemma \ref{a1}, we obtain that $P$ is a projective object.  \qed

\section{$n$-abelian quotient categories}
In this section, all categories are assumed to be $k$-liner Hom-finite, where $k$ is an algebraically closed field.
Let $\C$ be an additive category and $\X$ a subcategory of $\C$.
Recall that we say a morphism $f\colon A \to B$ in $\C$ is an $\X$-\emph{monic} if
$$\C(f,X)\colon \C(B,X) \to \C(A,X)$$
is an epimorphism for all $X\in\X$. We say that $f$ is an $\X$-\emph{epic} if
$$\C(X,f)\colon \C(X,A) \to \C(X,B)$$
is an epimorphism for all $X\in\X$.
Similarly,
we say that $f$ is a left $\X$-approximation of $B$ if $f$ is an $\X$-monic and $A\in\X$.
We say that $f$ is a right $\X$-approximation of $A$ if $f$ is an $\X$-epic and $B\in\X$.

We denote by $\C/\X$
the category whose objects are objects of $\C$ and whose morphisms are elements of
$\Hom_{\C}(A,B)/\X(A,B)$ for $A,B\in\C$, where $\X(A,B)$ is the subgroup of $\Hom_{\C}(A,B)$ consisting of morphisms
which factor through an object in $\X$.
Such category is called the (additive) quotient category
of $\C$ by $\X$. For any morphism $f\colon A\to B$ in $\C$, we denote by $\overline{f}$ the image of $f$ under
the natural quotient functor $\C\to\C/\X$.

\begin{remark}\label{rem}
If $\C$ has split idempotents, then $\C/\X$ has split idempotents, see the proof of \cite[Lemma 1.1(i)]{J}.
\end{remark}

\begin{lemma}\label{lem0}
Let $\C$ be an $n$-exangulated category and $\X$ an additive subcategory of $\C$.
If for any object $A\in\C$, there exists a distinguished $n$-exangle
$$A\xrightarrow{f_0}X_1\xrightarrow{f_1}X_2\xrightarrow{f_2}\cdots\xrightarrow{f_{n-1}}X_n\xrightarrow{f_n}B\overset{\delta}{\dashrightarrow}$$
where $f_0$ is a left $\X$-approximation of $A$ and $X_1,X_2,\cdots,X_n\in\X$. Then the quotient
category $\C/\X$ is a right $(n+2)$-angulated category with the following endofunctor and right $(n+2)$-angles:
\begin{itemize}
\item[\emph{(1)}] For any object $A\in\C$, we take a distinguished $n$-exangle
$$A\xrightarrow{~f_0~}X_1\xrightarrow{~f_1~}X_2\xrightarrow{~f_2~}\cdots\xrightarrow{~f_{n-1}~}X_n\xrightarrow{~f_n~} \mathbb{G}A\overset{\delta}{\dashrightarrow}$$
where $f_0$ is a left $\X$-approximation of $A$ and $X_1,X_2,\cdots,X_n\in\X$. Then $\mathbb{G}$ is a well-defined endofunctor of $\C/\X$.
\item[\emph{(2)}] For any distinguished $n$-exangle
 $$A_0\xrightarrow{g_0}A_1\xrightarrow{g_1}A_2\xrightarrow{g_2}A_3
 \xrightarrow{g_3}\cdots\xrightarrow{g_{n-1}}A_n\xrightarrow{g_n}A_{n+1}\overset{\eta}{\dashrightarrow}$$
where $g_0$ is an $\X$-monic, take the following commutative diagram of distinguished $n$-exangles
$$\xymatrix{
A_0 \ar[r]^{g_0}\ar@{=}[d]& A_1 \ar[r]^{g_1}\ar@{-->}[d]^{\varphi_1} & A_2 \ar[r]^{g_2}\ar@{-->}[d]^{\varphi_2}  & \cdots \ar[r]^{g_{n-1}}& A_n \ar[r]^{g_n}\ar@{-->}[d]^{\varphi_n}&A_{n+1}\ar@{-->}[r]^{\eta}\ar@{-->}[d]^{\varphi_{n+1}} & \\
A_0 \ar[r]^{f_0}&X_1 \ar[r]^{f_1} & X_2 \ar[r]^{f_2}  & \cdots \ar[r]^{f_{n-1}} & X_n \ar[r]^{f_n}&\mathbb{G}A_0\ar@{-->}[r]^{\delta}&.}$$
Then we have a complex
$$A_0\xrightarrow{~\overline{g_0}~} A_1\xrightarrow{~\overline{g_1}~}A_2\xrightarrow
{~\overline{g_2}~}\cdots\xrightarrow{~\overline{g_{n-1}}~}A_{n}\xrightarrow{~\overline{g_{n}}~}A_{n+1}\xrightarrow{~\overline{\varphi_{n+1}}~}\mathbb{G}A_0.$$ We define right $(n+2)$-angles in $\C/\X$ as the complexes which are isomorphic to complexes obtained in this way.
\end{itemize}
\end{lemma}

\proof  Since the proof is very similar to \cite[Theorem 3.11]{ZW}, we omit it.   \qed

\begin{lemma}\label{keylemma}
Let $\C$ be an $n$-exangulated category with split idempotents and $\X$ an additive subcategory of $\C$. Consider the following conditions:
\begin{itemize}
\item[\emph{(a)}] For any object $A\in\C$, there exist two distinguished $n$-exangles:
$$X_0\xrightarrow{~~}X_1\xrightarrow{~~}\cdots\xrightarrow{~~}X_{n-1}\xrightarrow{~~}X_n\xrightarrow{~f~}A\dashrightarrow$$
where $f$ is a right $\X$-approximation of $A$ and $X_0, X_1,\cdots,X_{n}\in\X$,
and
$$A\xrightarrow{~g~}X'_1\xrightarrow{~~}X'_2\xrightarrow{~~}\cdots\xrightarrow{~~}X'_n\xrightarrow{~~}X'_{n+1}\dashrightarrow$$
where $g$ is a left $\X$-approximation of $A$ and $X'_1,X'_2,\cdots,X'_{n+1}\in\X$.

\item[\emph{(b)}] For any distinguished $n$-exangle
 $$A_0\xrightarrow{f_0}A_1\xrightarrow{f_1}A_2\xrightarrow{f_2}A_3\xrightarrow{f_3}\cdots
 \xrightarrow{f_{n-1}}A_n\xrightarrow{f_n}A_{n+1}\overset{\delta}{\dashrightarrow}$$
in $\C$, $f_0$ is an $\X$-monic if $\overline{f_{n}}$ is an epimorphism in $\C/\X$ and
$f_{n}$ is an $\X$-epic if $\overline{f_{0}}$ is a monomorphism in $\C/\X$.
\end{itemize}
If the conditions {\rm (a)} and {\rm (b)} hold,  then $\C/\X$ is an $n$-abelian category.
\end{lemma}

\proof Since $\C$ has split idempotents, then $\C/\X$ has split idempotents by Remark \ref{rem}.
Thus (A0) is satisfied.

By Lemma \ref{lem0}, we obtain that
$(\C/\X,\mathbb{G})$ is a right $(n+2)$-angulated category.
By the construction of $\mathbb{G}$, we have $\mathbb{G}=0$.

For any morphism $\overline{f_0}\colon A_0\to A_1$, take a distinguished $n$-exangle
 $$A_0\xrightarrow{~g_0~}X'_1\xrightarrow{~g_1~}X'_2\xrightarrow{~g_2~}\cdots
 \xrightarrow{~g_{n-1}~}X'_n\xrightarrow{~g_n~}X'_{n+1}=\mathbb{G}A_0\overset{\delta}{\dashrightarrow}$$
with $g_0$ is a left $\X$-approximation of $A$ and $X'_1,X'_2,\cdots,X'_{n+1}\in\X$.
It is obvious that $\left[\begin{smallmatrix}
-g_0\\ f_0
\end{smallmatrix}
\right]\colon A_0\to X'_1\oplus A_1 $ is an $\X$-monic.
Note that $g_0$ is an inflation, by Lemma \ref{inflation}, we know that $\left[\begin{smallmatrix}
-g_0\\ f_0
\end{smallmatrix}
\right]$
is also an inflation.
So there exists a distinguished $n$-exangle
$$A_0\xrightarrow{\left[\begin{smallmatrix}
-g_0\\ f_0
\end{smallmatrix}
\right]}X'_1\oplus A_1 \xrightarrow{f_1}A_2\xrightarrow{f_2}A_3\xrightarrow{f_3}\cdots
\xrightarrow{f_{n-1}}A_n\xrightarrow{f_n}A_{n+1}\overset{\eta}{\dashrightarrow}.$$
By Lemma \ref{a1},  we get the following commutative diagram
$$\xymatrix{
A_0 \ar[r]^{\;\;\left[\begin{smallmatrix}
-g_0\\ f_0
\end{smallmatrix}
\right]\qquad}\ar@{=}[d]& X'_1\oplus A_1 \ar[r]^{\quad f_1}\ar[d]^{(-1,\thinspace 0)} & A_2 \ar[r]^{f_2}\ar@{-->}[d]^{\varphi_2}  & \cdots \ar[r]^{f_{n-1}}& A_n \ar[r]^{f_n}\ar@{-->}[d]^{\varphi_n}&A_{n+1}\ar@{-->}[r]^{\;\;\eta}\ar@{-->}[d]^{\varphi_{n+1}} & \\
A_0 \ar[r]^{g_0}&X'_1 \ar[r]^{g_1} & X'_2 \ar[r]^{g_2}  & \cdots \ar[r]^{g_{n-1}} & X'_n \ar[r]^{g_n}&X'_{n+1}\ar@{-->}[r]^{\;\;\delta}&
}$$
of distinguished $n$-exangles in $\C$.
It follows that $$A_0\xrightarrow{~\overline{f_0}~}A_1\xrightarrow{~\overline{f_1}~}A_2\xrightarrow{~\overline{f_2}~}\cdots\xrightarrow{~\overline{f_{n-2}}~}A_{n-1}
\xrightarrow{~\overline{f_{n-1}}~}A_{n}\xrightarrow{~\overline{f_{n}}~}A_{n+1}\xrightarrow{~~~}0$$
is a right $(n+2)$-angle in $\C/\X$.  By Lemma \ref{bu2}, we know that
$\overline{f_i}$ is a weak cokernel of $\overline{f_{i-1}},~i=1,2,\cdots,n$. Note that $\overline{f_n}$ is an epimorphism in $\C/\X$, we have that $\overline{f_n}$ is a cokernel of $\overline{f_{n-1}}$.  This shows that $(\overline{f_1},\overline{f_2},\cdots, \overline{f_{n+1}})$ is an $n$-cokernel of $\overline{f_0}$.

So any morphism in $\C/\X$ has $n$-cokernel. Dually we can show that $\C/\X$ has $n$-kernel.  Thus (A1) is satisfied.

If $\overline{f_0}$ is a monomorphism in $\C/\X$,
there exists a right $(n+2)$-angle
\begin{equation}\label{t00}
\begin{array}{l}
A_0\xrightarrow{~\overline{f_0}~}A_1\xrightarrow{~\overline{f_1}~}A_2
\xrightarrow{~\overline{f_2}~}\cdots\xrightarrow{~\overline{f_{n-2}}~}
A_{n-1}\xrightarrow{~\overline{f_{n-1}}~}A_{n}\xrightarrow{~\overline{f_{n}}~}
A_{n+1}\xrightarrow{~~~}0
\end{array}
\end{equation}
in $\C/\X$.
By the definition of the right $(n+2)$-angles in $\C/\X$,
without loss of generality, we may assume that the right $(n+2)$-angle (\ref{t00}) is induced by
the following commutative diagram
$$\xymatrix{
A_0 \ar[r]^{f_0}\ar@{=}[d]& A_1\ar[r]^{f_1}\ar[d]^{\varphi_1} & A_2 \ar[r]^{f_2}\ar[d]^{\varphi_2}  & \cdots \ar[r]^{f_{n-1}}& A_n \ar[r]^{f_n}\ar[d]^{\varphi_n}&A_{n+1}\ar@{-->}[r]^{\quad\eta}\ar[d]^{\varphi_{n+1}} &\\
A_0 \ar[r]^{g_0}&X_1 \ar[r]^{g_1} & X_2 \ar[r]^{g_2}  & \cdots \ar[r]^{g_{n-1}} & X_n \ar[r]^{g_n}&X_{n+1}\ar@{-->}[r]^{\quad\delta}&
}$$
of distinguished $n$-exangles, where $X_1,X_2,\cdots,X_{n+1}\in\X$.  Thus we have that
\begin{equation}\label{t01}
\begin{array}{l}
A_0\xrightarrow{f_0}A_1\xrightarrow{f_1}A_2\xrightarrow{f_2}A_3\xrightarrow{f_3}\cdots
\xrightarrow{f_{n-1}}A_n\xrightarrow{f_n}A_{n+1}\overset{\eta}{\dashrightarrow}.
\end{array}
\end{equation}
is a distinguished $n$-exangle in $\C$. By  hypothesis (b), we have that $f_{n}$ is an $\X$-epic since $\overline{f_0}$ is a monomorphism.
Thus the distinguished $n$-exangle (\ref{t01}) induces a left $(n+2)$-angle
$$0\xrightarrow{~~} A_0\xrightarrow{~\overline{f_0}~} A_1\xrightarrow{~\overline{f_1}~}A_2\xrightarrow
{~\overline{f_2}~}\cdots\xrightarrow{~\overline{f_{n-2}}~}A_{n-1}\xrightarrow{~\overline{f_{n-1}}~}A_{n}\xrightarrow{~\overline{f_{n}}~}A_{n+1}$$
in $\C/\X$.
By Lemma \ref{bu2} and its dual,  we have that
$$A_0\xrightarrow{~\overline{f_0}~}A_1\xrightarrow{~\overline{f_1}~}A_2\xrightarrow
{~\overline{f_2}~}\cdots\xrightarrow{~\overline{f_{n-2}}~}A_{n-1}\xrightarrow{~\overline{f_{n-1}}~}A_{n}\xrightarrow{~\overline{f_n}~}A_{n+1}$$
is an $n$-exact sequence in $\C/\X$.  Thus (A2) is satisfied. Dually we can show that (A2)$^{\textrm{op}}$ is satisfied.

Therefore $\C/\X$ is an $n$-abelian category. \qed
\medskip

Motivated by the definition of cluster tilting subcategories of $(n+2)$-angulated categories,
which is due to Zhou-Zhu \cite[Definition 3.3]{ZZ3}.
We define cluster tilting subcategories of $n$-exangulated categories.
\begin{definition}\label{defn}
Let $\C$ be an $n$-exangulated category and $\X$ an additive subcategory of $\C$.
$\X$ is called \emph{cluster-tilting} if

(1) $\E(\X,\X)=0$.

(2)  For any object $C\in\C$, there are two  distinguished $n$-exangles
$$X_0\xrightarrow{~~}X_1\xrightarrow{~~}\cdots\xrightarrow{~~}X_{n-1}\xrightarrow{~~}X_n\xrightarrow{~~}C\dashrightarrow$$
where $X_0, X_1,\cdots,X_{n}\in\X$ and
$$C\xrightarrow{~~}X'_1\xrightarrow{~~}X'_2\xrightarrow{~~}\cdots\xrightarrow{~~}X'_n\xrightarrow{~~}X'_{n+1}\dashrightarrow$$
where $X'_1,X'_2,\cdots,X'_{n+1}\in\X$.\end{definition}

\begin{remark}
When $\C$ is an $(n+2)$-angulated category, this definition is just the Definition 3.3 in \cite{ZZ3}.
\end{remark}

\begin{remark}
Let $\C$ be an $n$-exangulated category with enough projectives and enough injectives
and $\X$ a cluster tilting subcategory of $\C$. By Lemma \ref{projective} and its dual,
we know that $\P\subseteq\X$ and $\I\subseteq\X$.
\end{remark}

Our main result is the following.

\begin{theorem}\label{main}
Let $\C$ be an $n$-exangulated category with split idempotents and $\X$ a cluster-tilting subcategory of $\C$.
 If $\C$ has enough projectives and enough injectives,
then $\C/\X$ is an $n$-abelian category.
\end{theorem}

\proof
We only need to show that cluster tilting subcategories satisfy the conditions (a) and (b) in Lemma \ref{keylemma}.

By the definition of cluster-tilting subcategories, the condition (a) holds in Lemma \ref{keylemma}.

For any  distinguished $n$-exangle
 $$A_0\xrightarrow{f_0}A_1\xrightarrow{f_1}A_2\xrightarrow{f_2}\cdots\xrightarrow{f_{n-1}}A_n\xrightarrow{f_n}
 A_{n+1}\overset{\eta}{\dashrightarrow}$$
in $\C$.

Now assume that $\overline{f_{n}}$ is an epimorphism in $\C/\X$.
Since $\C$ has enough injectives, there exists a distinguished $n$-exangle
 $$A_0\xrightarrow{g_0}I_1\xrightarrow{g_1}I_2\xrightarrow{g_2}\cdots\xrightarrow{g_{n-1}}I_n\xrightarrow{g_n}B
 \overset{\delta}{\dashrightarrow}$$
where $I_1,I_2,\cdots,I_n\in\I$. Since $I_1$ is an injective object, there exists
a morphism $\varphi_1\colon A_1\to I_1$ such that $\varphi_1f_0=g_0$.
By Lemma \ref{a1},  we get the following commutative diagram
\begin{equation}\label{tt3}
\begin{array}{l}
\xymatrix{
A_0 \ar[r]^{f_0}\ar@{=}[d]& A_1\ar[r]^{f_1}\ar[d]^{\varphi_1} & A_2 \ar[r]^{f_2}\ar@{-->}[d]^{\varphi_2}
& \cdots \ar[r]^{f_{n-1}}& A_n \ar[r]^{f_n}\ar@{-->}[d]^{\varphi_n}&A_{n+1}\ar@{-->}[r]^{\quad\eta}\ar@{-->}[d]^{\varphi_{n+1}} &\\
A_0 \ar[r]^{g_0}&I_1 \ar[r]^{g_1} & I_2 \ar[r]^{g_2}  & \cdots \ar[r]^{g_{n-1}} & I_n \ar[r]^{g_n}&B\ar@{-->}[r]^{\quad\delta}&
}
\end{array}
\end{equation}
of distinguished $n$-exangles.
It follows that $\varphi_{n+1}f_n=g_n\varphi_n$.
Since $I_n\in\I\subseteq\X$, we have $\overline{\varphi_{n+1}}\circ\overline{f_n}=0$ in $\C/\X$.
Since $\overline{f_{n}}$ is an epimorphism, we get
$\overline{\varphi_{n+1}}=0$. That is to say,
there exist morphisms
$u\colon A_{n+1}\to X_0 $ and $v\colon X_0\to B$ where $X_0\in\X$ such that $\varphi_{n+1}=vu$.

Now we prove that $f_0$ is an $\X$-monic. Let $a\colon A_0\to X$ be any morphism in $\C$ where $X\in\X$.
Since $\C$ has enough injectives, there exists a distinguished $n$-exangle
 \begin{equation}\label{tt4}
\begin{array}{l}
  X\xrightarrow{h_0}I'_1\xrightarrow{h_1}I'_2\xrightarrow{h_2}\cdots\xrightarrow{h_{n-1}}I'_n\xrightarrow{h_n}C
 \overset{\theta}{\dashrightarrow}
 \end{array}
\end{equation}
where $I'_1,I'_2,\cdots,I'_n\in\I$. Since $I_1$ is an injective object, there exists
a morphism $\psi_1\colon I_1\to I'_1$ such that $\psi_1g_0=h_0a$.
By Lemma \ref{a1},  we get the following commutative diagram
\begin{equation}\label{tt5}
\begin{array}{l}
\xymatrix{
A_0 \ar[r]^{g_0}\ar[d]^a& I_1\ar[r]^{g_1}\ar[d]^{\psi_1} & I_2 \ar[r]^{g_2}\ar@{-->}[d]^{\psi_2}
& \cdots \ar[r]^{g_{n-1}}& I_n \ar[r]^{g_n}\ar@{-->}[d]^{\psi_n}&B\ar@{-->}[r]^{\quad\delta}\ar@{-->}[d]^{\psi_{n+1}} &\\
X \ar[r]^{h_0}&I'_1 \ar[r]^{h_1} & I'_2 \ar[r]^{h_2}  & \cdots \ar[r]^{h_{n-1}} & I'_n \ar[r]^{h_n}&C\ar@{-->}[r]^{\quad\theta}&
}
\end{array}
\end{equation}
of distinguished $n$-exangles.
By composing the commutative
diagrams (\ref{tt3}) and (\ref{tt5}), we have the commutative diagram
\begin{equation}\label{tt6}
\begin{array}{l}
\xymatrix{
A_0 \ar[r]^{f_0}\ar[d]^a& A_1\ar[r]^{f_1}\ar[d]^{\psi_1\varphi_1} & A_2 \ar[r]^{f_2}\ar[d]^{\psi_2\varphi_2}
& \cdots \ar[r]^{f_{n-1}}& A_n \ar[r]^{f_n}\ar[d]^{\psi_n\varphi_n}&A_{n+1}\ar@{-->}[r]^{\quad\eta}\ar[d]^{\psi_{n+1}\varphi_{n+1}} &\\
X \ar[r]^{h_0}&I'_1 \ar[r]^{h_1} & I'_2 \ar[r]^{h_2}  & \cdots \ar[r]^{h_{n-1}} & I'_n \ar[r]^{h_n}&C\ar@{-->}[r]^{\quad\theta}&
}
\end{array}
\end{equation}
of distinguished $n$-exangles.
Applying the functor $\C(X_0,-)$ to the distinguished $n$-exangle (\ref{tt4}), we have the following
exact sequence:
$$\C(X_0,I'_n)\xrightarrow{\C(X_0,h_n)}\C(X_0,C)\xrightarrow{~}\E(X_0,X).$$
Since $\X$ is cluster tilting, we have $\E(X_0,X)=0$.
This shows that $\C(X_0,h_n)$ is an epimorphism.
So there exists a morphism $w\colon X_0\to I'_n$ such that $\psi_{n+1} v=h_nw$.
It follows that $\psi_{n+1}\varphi_{n+1}=h_n(w\varphi_{n+1})$.
Namely, $\psi_{n+1}\varphi_{n+1}$ factors through $h_n$.
By Lemma \ref{factor}, we get that $a$ factors through $f_0$, that is,
there exists a morphism $b\colon A_1\to X$ such that $a=bf_0$.
This shows that $f_0$ is an $\X$-monic.

Dually we can show that if $\overline{f_{0}}$ is a monomorphism in $\C/\X$,
then $f_{n}$ is an $\X$-epic.

Therefore the condition (b) holds in Lemma \ref{keylemma}.
By Lemma \ref{keylemma}, we have that $\C/\X$ is an $n$-abelian category.  \qed

\medskip
\begin{remark}
An important source of examples of $n$-abelian categories are $n$-cluster-tilting subcategories.
That is, Jasso \cite[Theorem 3.6]{Ja} showed that an $n$-cluster-tilting subcategory of
an abelian category are an $n$-abelian category.  Later, Zhou-Zhu \cite[Theorem 3.4]{ZZ3} proved that
any quotient of an $(n+2)$-angulated category modulo a cluster tilting subcategory is an $n$-abelian category.
Our main result gives a new source to construct $n$-abelian categories.
\end{remark}

\begin{remark}
In Theorem \ref{main}, when $n=1$, it is just the Theorem 3.4 in \cite{ZZ2}.
\end{remark}

Applying Theorem \ref{main} to $(n+2)$-angulated categories, we obtain
the following result.

\begin{corollary}\label{cor1}{\rm\cite[Theorem 3.4]{ZZ3}}
Let $\C$ be an $(n+2)$-angulated category with split idempotents and $\X$ a cluster tilting subcategory of $\C$.
Then $\C/\X$ is an $n$-abelian category.
\end{corollary}

\proof Since any $(n+2)$-angulated category can be viewed as an $n$-exangulated category with
enough projectives and enough injectives. This follows from  Theorem \ref{main}.  \qed

\begin{remark}
In Corollary \ref{cor1}, when $n=1$, it is just the Theorem 3.3 in \cite{KZ}.
\end{remark}

Applying Theorem \ref{main} to $n$-exact categories, we obtain
the following result.

\begin{corollary}\label{cor2}
Let $\C$ be an $n$-exact category with split idempotents and $\X$ a cluster-tilting subcategory of $\C$.
If $\C$ has enough projectives and enough injectives,
then $\C/\X$ is an $n$-abelian category.
\end{corollary}

\proof Since any $n$-exact category can be viewed as an $n$-exangulated category.
This follows from  Theorem \ref{main}.  \qed

\begin{remark}
In Corollary \ref{cor2}, when $n=1$, it is just the Theorem 3.3 in \cite{DL}.
\end{remark}

\textbf{Yu Liu}\\
School of Mathematics, Southwest Jiaotong University, 610031, Chengdu,
Sichuan, P. R. China\\
E-mail: \textsf{liuyu86@swjtu.edu.cn}\\[0.3cm]
\textbf{Panyue Zhou}\\
College of Mathematics, Hunan Institute of Science and Technology, 414006, Yueyang, Hunan, P. R. China.\\
E-mail: \textsf{panyuezhou@163.com}

\end{document}